\documentclass[11pt]{amsart}
\usepackage{amssymb,latexsym,comment,url,amsmath}

\usepackage[T1]{fontenc}
\usepackage{rotating,graphicx}
\usepackage{xcolor}
\usepackage{diagbox} 
\usepackage{mathtools}

\usepackage{longtable}

\newcommand{\Char}{\operatorname{char}}

\newcommand{\Q}{{\mathbb Q}}

\newcommand{\Z}{{\mathbb Z}}

\newcommand{\Magma}{{\sf Magma}}

\newcommand{\Mathematica}{{\sf Mathematica}}

\newenvironment{Proof}{\par\noindent{\sc Proof:}}%
                      {\hspace*{\fill}\nobreak$\Box$\par\medskip}
                       {\hspace*{\fill}\nobreak$\Box$\par\medskip}

\newtheorem{Proposition}{Proposition}[section]
\newtheorem{Theorem}[Proposition]{Theorem}
\newtheorem{Lemma}[Proposition]{Lemma}

\newtheorem{Corollary}[Proposition]{Corollary}

\theoremstyle{definition}
\newtheorem{Definition}[Proposition]{Definition}
\newtheorem{Remark}[Proposition]{Remark}
 \newtheorem{Example}[Proposition]{Example}

\renewcommand{\baselinestretch}{1.1}

\begin{document}

\title[Rational torsion on hyperelliptic jacobian varieties]%
{Rational torsion on hyperelliptic jacobian varieties}

\author[M. Sadek]%
{Mohammad~Sadek}
\email{mohammad.sadek@sabanciuniv.edu}
\address{Faculty of Engineering and Natural Sciences, Sabanc{\i} University, Tuzla, \.{I}stanbul, 34956 Turkey}

\author[H. Suluyer]%
{Hamide~Suluyer}
\email{hamide.kuru@atilim.edu.tr, hamidekuru@sabanciuniv.edu}

\address{Department of Mathematics, At{\i}l{\i}m University, 06830 G\"olba\c{s}{\i}, Ankara, Turkey}
\address{Faculty of Engineering and Natural Sciences, Sabanc{\i} University, Tuzla, \.{I}stanbul, 34956 Turkey}

\begin{abstract}
{ It was conjectured
by Flynn that there exists a constant $\kappa$ such that, for any integer $g \ge 2$, any $m \le \kappa g$, there
exists a hyperelliptic curve of genus $g$ over $\Q$ with a rational $m$-torsion point on its Jacobian. Lepr\'{e}vost proved this conjecture with $\kappa=3$. In this work we prove that given an integer $N$ in the interval $[3g,4g+1]$, $g\ge 3$, satisfying certain partition conditions, there exist parametric families of hyperelliptic Jacobian varieties with a rational torsion point of order $N$. In particular, we establish the existence of such varieties for $N=4g+1$ when $g$ is odd and for $N=4g-1$ when $g$ is even. A few explicit applications of this result produce the first known infinite examples of torsion $13$ when $g=3$, torsion $15$ when $g=4$, and torsion $17,18,21$ when $g=5$. In fact, we show that infinitely many of the latter abelian varieties are absolutely simple. 
}
\end{abstract}

\maketitle

\let\thefootnote\relax\footnotetext{ \textbf{Keywords:} Hyperelliptic curves, Jacobian varieties, torsion\\

\textbf{2020 Mathematics Subject Classification:} 11G30, 14H25}

\section{Introduction}
Given an integer $n\ge 1$, the existence of a uniform bound, that depends solely on $n$, on the size of the torsion subgroups of elliptic curves defined over a number field of degree $n$ is guaranteed by a result of Merel, \cite{Merel}. Complete lists for the structure of these torsion subgroups are explicitly given when $n=1,2,3$ in \cite{Maz78,Kam92, KM,D}. For abelian varieties of fixed dimension $d\ge 2$, not only such a uniform boundedness result is far from reach, but no complete classification of possible torsion orders has been realized over any number field.   

A computational approach toward studying torsion subgroups of abelian varieties of dimension $d\ge 2$ is to find integers that are realized as orders of torsion points on these varieties. For example, a family of genus $2$ curves whose Jacobian varieties possess torsion points of order $\ell$ over $\Q$ was given in \cite{Flynn1} when $\ell=11$, and in \cite{Lep13,Lep15,Lep2}, when $\ell= 13$, $15$, $17$, $19$ or $21$. 
In \cite{Howe}, several explicit families of curves of genus $2$ and $3$ whose Jacobians possess large rational torsion subgroups were displayed.

Researchers questioned the possibility of constructing families of abelian varieties in which the size of the rational torsion subgroup
of the variety increases as the dimension of the variety increases. In \cite{Flynn2}, fixing an even integer $g\ge 2$, it was proved that any integer in the interval $[g^2+2g+1,g^2+3g+1]$ is realized as the order of a $\Q$-rational torsion point on the Jacobian of a hyperelliptic curve of genus $g$. Other families of hyperelliptic Jacobian varieties of dimension $g$ with $\Q$-rational torsion points of quadratic order in $g$ were given in \cite{KS, Lep1,Lep3}. In \cite{Flynn1}, it was conjectured that there exists a constant $\kappa$, independent of $g$, such that for every $m \le \kappa g$, there exists an infinite family of hyperelliptic curves of genus $g$ over $\Q$ with a torsion point of order $m$ on its Jacobian. Lepr\'{e}vost proved this conjecture with $\kappa=3$ in \cite{LepCon}.
One may ask whether the integers between the
linear bounds introduced by Lepr\'{e}vost and the quadratic bounds mentioned earlier are realizable as torsion orders of rational points of hyperelliptic Jacobian varieties. This manuscript attempts to address
this gap in bounds and introduce new larger linear bounds.

Since $11$ is the smallest prime for which there is no
elliptic curve over $\Q$ with rational point of order the prime, finding abelian varities with torsion order $11$ was an interesting challenge for researchers. Jacobians of genus $2$ curves with torsion order $11$ over $\Q$ can be found in \cite{Flynn1,Bernard,Dthesis,DD}; Flynn's conjecture, later proved by Lepr\'{e}vost, establishes the existence of infinitely many such varieties when $g\ge 4$. Recently, Jacobians of genus $3$ hyperelliptic curves with a $\Q$-rational torsion point of order $11$ was given in \cite{Dthesis,D}. This makes $13$ the next challenging prime for abelian varieties of dimension $g$ when $g=3,4$.  

In this work, we extend the interval given by Lepr\'{e}vost to prove Flynn's conjecture to include new integers that are realized as torsion orders for infinitely many hyperelliptic Jacobian varieties.  The
divisor at infinity has finite order in the Jacobian of a hyperelliptic curve if and only
if a certain continued fraction expansion of a power series is periodic. We show that for every integer $N$ in the interval $[3g,4g+1]$, $g\ge 3$, satisfying certain partition conditions, there exist parametric families of hyperelliptic curves of genus $g$ whose Jacobian varieties have a $\Q$-rational torsion point of order $N$. These partition conditions appear due to the fact that the order of the divisor at infinity in the Jacobian
is the sum of the degrees of some partial quotients in the corresponding continued fraction expansion. Some of the integers satisfying these partition conditions are $4g+1$ when $g$ is odd; and $4g-1$ when $g$ is even. Consequently, we produce infinite families of hyperelliptic Jacobian varieties with torsion order $N=13$ when $g=3$, torsion $15$ when $g=4$, and torsion $17,18,21$ when $g=5$. These are the first examples of infinitely many varieties corresponding to the aforementioned pairs $(g,N)$. For a survey of all torsion points of a certain order that have been found for Jacobian varieties of
hyperelliptic curves of genus $2, 3$, and $4$, the reader may consult \cite{Cthesis}. In addition, we show that infinitely many of the Jacobian varieties that we construct are absolutely simple. 

The idea that orders of torsion points on abelian varieties are linked to periodic continued fractions goes back to Abel and Chebychev. The work of Adams and Razar, \cite{Adams}, explained this link algebraically with much emphasis on torsion points of elliptic curves. However, continued fractions have not been used much to study torsion points on abelian varieties of higher dimension. We show that continued fractions over function fields can be employed to obtain new results on realizable torsion orders on Jacobian varieties of a fixed dimension.

\subsection*{Acknowledgment}
The authors would like to thank the anonymous referees for the thorough reading of the manuscript and for many useful suggestions and comments  that improved the manuscript.
 All the calculations in this work were performed using \Magma, \cite{Magma}, and \Mathematica, \cite{Mathematica}. 
 This work builds upon the findings presented in the doctoral dissertation of Hamide Suluyer, \cite{Suluyer}.
  M. Sadek is supported by The Scientific and Technological Research Council of Turkey, T\"{U}B\.{I}TAK, research grant ARDEB 1001/122F312.

\section{Torsion divisors at infinity}
Throughout this section $k$ is a field with $\Char k\ne 2$. We let $\overline{k}$ be an algebraic closure of $k$. The material in this section can be found in \cite{Adams}. Let \( C \) be a curve defined over a field \( k \).
\subsection{Continued fractions}

   Let \( O \in C(k) \) be a non-singular \( k \)-rational point on \( C \) and \( K = k(C) \) be the corresponding function field. The valuation corresponding to $O$ will be denoted by $\operatorname{ord}_O$.
  Let \( x \in K \) be a function with a simple pole at \( O \). Then \( x^{-1} \) is a uniformizing parameter at \( O \). For any \( \alpha \in K \), there is a unique $\beta\in k[x]$ such that \( \operatorname{ord}_O(\alpha-\beta)>0 \). We denote $\beta=[\alpha]$. We notice that the completion of $K$ at $O$ is the power series
field $k((x^{-1}))$. If $\alpha$ is expanded as a Laurent series in $x^{-1}$, then $[\alpha]$ is
the non-positive part of the expansion, i.e., $[\alpha]$ is a polynomial in $x$. 
  
Define a function
$$\begin{aligned} \varphi: K-k[x]  \rightarrow K, \qquad \varphi(\alpha) =1 /(\alpha-[\alpha]) \end{aligned}.$$ For positive integers $r$ let $\varphi_{r}$ be the $r$-fold composition of $\varphi$ with itself and let $\varphi_{0}$ be the identity map, namely,
$$
\alpha_{0}=\alpha,\quad \alpha_{r} =\varphi_{r}(\alpha),
  \quad \alpha_{r+1} =1 /\left(\alpha_{r}-\left[\alpha_{r}\right]\right)\quad \text { for } r \geqslant 1 .$$ 
  In addition,
$a_{r} :=\left[\alpha_{r}\right]$ is called the {\em partial quotient of $\alpha$}.
Now the continued fraction expansion of $\alpha$ is given by $$
\alpha=a_0 + \cfrac{1}{a_1
          + \cfrac{1}{a_2 + {}_{\ddots}} }.
$$
One can write the latter continued fraction in short as 
$\alpha=\left[a_{0}; a_{1}, a_{2}, \ldots\right]$.

\begin{Definition}
   Let $a_i$ be the $i$-th partial quotient of the continued fraction expansion of $y\in K$.
The continued fraction of $y$ is said to be periodic if there is a positive integer $n$ such that $a_{h+n}=a_h$ for every $h\ge 1$. The smallest such $n$ is called the period of the continued fraction. 
\end{Definition}

\begin{Definition}
A sequence $(a_1,\cdots,a_{m-1})$ is said to be {\em skew symmetric} of skew value $\gamma\in k-\{0\}$ if $a_{(m-i)}=c_i a_{i}$, $ 1\le i< (m+1)/2$, where $c_i=\gamma$ if $i$ is even; and $c_i=1/\gamma$ otherwise. 
\end{Definition}

\subsection{The divisor at infinity}
Let $C$ be a hyperelliptic curve of genus $g\ge 2$ over $k$ given by the equation $y^2=f(x)$ where $f(x)\in k[x]$ is a separable polynomial of even degree $2g+2$, and the leading coefficient of $f(x)$ is a nonzero square in $k$. 
Let $\infty_+$ and $\infty_-$ be the two $k$-rational points at infinity on $C$. 

There exists a natural embedding of the curve $C$ into its Jacobian $J$ that maps a point $P$ to the divisor class $[P - D]$, where $D$ is a fixed divisor of degree $1$. This map restricts to $C(k)\xhookrightarrow{} J(k)$ if $D$ is a $k$-rational divisor. In particular, if there exists a point $P_0 \in C(k)$, one can choose $D$ to be the rational divisor $P_0$.

The {\em divisor at infinity} on $C$ is the divisor $D_{\infty}=\infty_+-\infty_-$. It defines a point in $J(k)$. We say that $D_{\infty}$ is {\em torsion} of order $N$ if its class in $J$ has order $N$. The following result can be found in \cite[Theorem 5.1]{Adams} when $g=1$, and \cite[\S 4 and \S 5]{Poortentran} when $g>1$.
\begin{Proposition}\label{prop:1}
Let $C$ be a hyperelliptic curve described as above.
 The divisor $D_{\infty}$ is torsion if and only if 
the continued fraction of $\sqrt{f(x)}$ is periodic.

In the latter case, the following statements hold:
\begin{itemize}
\item[a)] $\sqrt{f(x)}=[a_0;\overline{a_1,a_2,\cdots,a_{m-1},2\gamma a_0,a_{m-1},\cdots, a_2,a_1,2a_0}]$, where $(a_1,a_2,\cdots,a_{m-1})$ is skew symmetric of skew value $\gamma^{-1}$. If $m$ is even, then $\gamma=1$, hence the period of the continued fraction is $m$.
\item[b)] $\deg a_0=g+1$, and $1\le \deg a_i\le g$ for $1\le i\le m-1$.
\item[c)] The order of $D_{\infty}$ is $g+1+\sum_{i=1}^{m-1}\deg a_i$.
\end{itemize}
\end{Proposition}

In Proposition \ref{prop:1}, we remark that the fact that $D_{\infty}$ is torsion gives rise to a unit $p(x)+q(x)\sqrt{f(x)}$ in $k(x,\sqrt{f(x)})$ where $p(x)/q(x)=[a_0;a_1,a_2,\cdots,a_{m-1}]$. Moreover, the order of $D_{\infty}$ is the degree of $p(x)$. The interested reader may consult \cite{Poortentran} for more details.

\begin{Remark}
\label{Rem:1}
In Proposition \ref{prop:1}, we call (the minimal) $m$ the {\em quasi-period length}. Recalling that $a_{(m-i)}= a_{i}$, $ 1\le i< (m+1)/2$, if $\gamma=1$, then the continued fraction expansion of $\sqrt{f(x)}$ is $[a_0;\overline{a_1,a_2,\cdots,a_{m-1},2 a_0,a_{1},\cdots,a_{m-1},2a_0}]$. One sees now that for $\gamma=1$, the (minimal) period length $n$ is equal to the quasi-period length $m$. Further, if $m\ne n$, then $n=2m$ and $m$ must be odd. In the latter case, $m$ is called a {\em strict} quasi-period length.
\end{Remark}

\section{Construction of curves}
Fixing an integer $r\ge 1$, we set $\Q(\textbf{w}):=\Q(w_1,\cdots,w_r)$, where $w_1,\cdots,w_r$ are indeterminates. In this section, we construct hyperelliptic curves $C$ of genus $g$ over $\Q(\textbf{w})$ described by the equation
\begin{eqnarray*}
y^2=h_{2g+2}(\textbf{w})^2x^{2g+2}+h_{2g+1}(\textbf{w})x^{2g+1}+\cdots+h_{0}(\textbf{w}),\qquad h_{i}(\textbf{w})\in\Q(\textbf{w}),\quad i=0,\cdots,2g+2,\end{eqnarray*}
where the discriminant of the latter equation is nonzero.
The curve $C$ admits two $\Q(\textbf{w})$-rational points at infinity denoted by $\infty^+$ and $\infty^-$. 

We investigate the case when the divisor $D_{\infty}=\infty^+-\infty^-$ in the Jacobian of $C$ is torsion with the quasi-period length $m$ and the period length $n$ are both equal to $6$. In view of Proposition \ref{prop:1}, this implies that the continued fraction expansion of $y $ is given by 
\begin{eqnarray}
\label{eq1}
[a_0(x);\overline{a_1(x), a_2(x), a_3(x), a_2(x), a_1(x),2a_0(x)}]
\end{eqnarray}
where $\deg a_0=g+1$ and $1\le \deg a_i\le g$ when $i=1,2,3$. 

The reason why we consider the case $m=n=6$ is that our attempts to construct hyperelliptic curves for which the corresponding quasi-period length $m\le 5$ produced hyperelliptic curves whose Jacobians possess torsion divisors at infinity with order $O(kg)$, $k < 3$, see \cite[Chapter 3]{Suluyer}. We recall that Lepr\'{e}vost, \cite{LepCon}, established the existence of hyperelliptic curves of genus $g$ whose Jacobians possess rational torsion points of order $3g$.

The following lemma discusses the existence of polynomial solutions to a certain Diophantine equation. This equation will appear during the course of the proof of the main theorem of this section, Theorem \ref{thm1}, in which we construct our desired curves. 

\begin{Lemma}
\label{Lem}
Let \( r , a_1   \in \mathbb{Q}[x] \). Consider the following polynomial Diophantine equation
\[
-2q -2a_2 +a_1 r +a_3 r -2q a_1 a_2 +a_1 a_2 a_3  r -a_2^2a_3 = 0
.\]
Then
\[
\begin{aligned}
q  &= \frac{u r }{2}, \\
a_2  &= a_1 r , \\
a_3  &= a_1 ^2 u r  + a_1  + u,
\end{aligned}
\]
where $u\in\Q[x]$, provide a solution to this Diophantine equation. 
\end{Lemma}
\begin{Proof}
One may check that the claimed solution satisfies the polynomial Diophantine equation.
\end{Proof}
\begin{Theorem}
\label{thm1}
Fix an integer $g\ge 3$. Let $\alpha,\beta\ge 1,\gamma\ge 0$ be integers such that $2\alpha+2\beta+\gamma=g+1$. Let $a_1 ,r ,u \in \Q[x]$ be of degrees $\alpha,\beta,\gamma$, respectively. If the affine equation 
 \[y^2=r ^2\left(u (a_1 ^2r +1)+a_1 \right)^2+4\left(u a_1 r^2 +r \right)\]
describes a hyperelliptic curve $C$, then the divisor at infinity is torsion of order $g+1+6\alpha+3\beta+\gamma$. 

In addition, if $\beta=1$, and the Galois group of the polynomial $$r \left(u (a_1 ^2r +1)+a_1 \right)^2+4\left(u a_1 r +1\right)$$ is either the full symmetric group $\operatorname{\textbf{S}}_{2g+1}$ or the alternating group $\operatorname{\textbf{A}}_{2g+1}$, then the Jacobian $J_C$ of $C$ satisfies that $\operatorname{End}(J_C) = \Z$, in
particular, $J_C$ is an absolutely simple abelian variety.
\end{Theorem}
\begin{Proof}
In view of (\ref{eq1}), a hyperelliptic curve \( C \) for which the divisor at infinity is torsion with \( m = n = 6 \) is described. Expanding the continued fraction (\ref{eq1}) as a rational function yields that $C$ is described by an equation of the form

\[
y^2 = a_0^2 + \frac{a_2 (2 + a_2 a_3) + 
 2 a_0 (1 + a_2 a_3 + a_1 a_2 (2 + a_2 a_3))}{(1 + a_1 a_2)(a_3 + 
   a_1 (2 + a_2 a_3))}.
\]
   In order for the latter equation to be a polynomial equation, the latter fraction must be a polynomial $h\in \Q[x]$, i.e., one has $$h=\frac{a_2 (2 + a_2 a_3) + 
 2 a_0 (1 + a_2 a_3 + a_1 a_2 (2 + a_2 a_3))}{(1 + a_1 a_2)(a_3 + 
   a_1 (2 + a_2 a_3))}\in\Q[x].$$ This yields that $$a_0=\frac{a_1 h}{2}+\frac{-a_2^2 a_3 + (a_1 + a_3) h + 
 a_2 (-2 + a_1 a_3 h)}{2 (1 + a_2 a_3 + a_1 a_2 (2 + a_2 a_3))}.$$ 
 We remark that if $\deg h\ge\deg a_2$, then the highest degree term of the numerator of the latter fraction is $a_1a_2a_3h$. In addition, the degree of the numerator will be larger than the degree of the denominator; which allows further division. Therefore, we assume that $h=2q a_2+r$ with $\deg r<\deg a_2$. We then obtain that
 \[a_0=\frac{a_1 h +2q}{2}+\frac{-2a_2-2q+a_1r+a_3r-2qa_1a_2+a_1a_2a_3 r-a_2^2a_3}{2 (1 + a_2 a_3 + a_1 a_2 (2 + a_2 a_3))}.\] We notice that if the degree of $q $ is strictly less than the sum of the degrees of $a_2 $ and $a_3 $, then the numerator of the latter fraction is of strictly less degree than the denominator. Therefore, for $a_0$ to be a polynomial in $\Q[x]$, we may assume that  $-2a_2-2q+a_1r+a_3r-2qa_1a_2+a_1a_2a_3 r-a_2^2a_3= 0\in \Q[x]$. According to Lemma \ref{Lem}, a solution to this polynomial Diophantine equation is provided giving rise to the affine equation $y^2=(a_1h+2q)^2/4 +h$. After using the explicit description of the solution in Lemma \ref{Lem}, the curve $C$ is described by
 \[y^2=r ^2\left(u (a_1 ^2r +1)+a_1 \right)^2+4\left(u a_1 r^2 +r \right).\]
 Now, when $\beta=1$, the rest of the statement will follow from \cite[Theorem 1.3]{Zarhin}.
\end{Proof}
\begin{Remark}
In Theorem \ref{thm1}, we insist that $\beta:=\deg r\ge 1$, since otherwise $a_3=2a_0$. This means that in (\ref{eq1}), the quasi period length, $m$, and the period length, $n$, satisfy that $m=n=3$. 
\end{Remark}
We notice that in Theorem \ref{thm1}, if $\displaystyle a_1 =\sum_{i=0}^{\alpha} A_ix^i$, $\displaystyle r =\sum_{i=0}^{\beta}R_ix^i$, and $\displaystyle u =\sum_{i=0}^{\gamma} U_ix^i$, then the curve $C$ is defined by an equation of the form $y^2= f(x) $ where $$f(x)\in \Q[A_0,\cdots A_{\alpha}, R_0,\cdots R_{\beta}, U_0,\cdots,U_{\gamma}][x].$$ 
The discriminant of the equation $y^2=f(x)$ describing $C$ is a constant multiple of the discriminant $\Delta_f$ of the polynomial $f$, see for example \cite{Loc} or \cite{Liu}, hence it is a polynomial in $\Q[A_0,\cdots A_{\alpha}, R_0,\cdots R_{\beta}, U_0,\cdots,U_{\gamma}]$. 
If we consider the affine space $\mathbb{A}^{\alpha+\beta+\gamma+3}$ with coordinates $A_0,\cdots A_{\alpha}$, $R_0,\cdots R_{\beta}$, $U_0,\cdots,U_{\gamma}$, then for a choice 
 $(A_0,\cdots A_{\alpha}, R_0,\cdots R_{\beta}, U_0,\cdots,U_{\gamma})\in \mathbb{A}^{\alpha+\beta+\gamma+3}$, the equation $y^2=f(x)$ defines a hyperelliptic curve of genus $g$ over $\Q$ if and only if $(A_0,\cdots A_{\alpha}, R_0,\cdots R_{\beta}, U_0,\cdots,U_{\gamma})$ is chosen to lie in the complement of the hypersurface $\Delta_f=0$ in the affine space $\mathbb{A}^{\alpha+\beta+\gamma+3}$.

As a direct consequence of Theorem \ref{thm1} and the paragraph above, we get the following result. 
\begin{Corollary}
\label{cor1}
Fix an integer $g\ge 3$. Let $N$ be an integer such that $3g\le N\le 4g+1$. If there exists integers $\alpha,\beta\ge 1,\gamma\ge 0$ such that $2\alpha+2\beta+\gamma=g+1$ and $N=g+1+6\alpha+3\beta+\gamma$, then there exists an $(\alpha+\beta+\gamma+3)$-parameter family of  hyperelliptic curves of genus $g$ over $\Q$ whose Jacobians possess a rational point of order $N$. In particular, if $g$ is odd, respectively even, there is an $(\alpha+\beta+\gamma+3)$-parameter family of hyperelliptic curves of genus $g$ over $\Q$ whose Jacobians possess a rational point of order $4g+1$, respectively $4g-1$.
\end{Corollary}

 Although our trials to study hyperelliptic curves $C$ for which the corresponding continued fraction in Proposition \ref{prop:1} has quasi-period length $m\le 5$ yielded torsion divisors at infinity with order $\operatorname{O}(kg)$, where $k<3$, we successfully obtained families of hyperelliptic curves whose jacobain contains torsion points whose order is not listed in literature. As an example, when $m=n=4$, one has $y=[a_0;\overline{a_1,a_2,a_1,2a_0}]$. In other words, one sees that $$y^2=a_0^2+\frac{2a_0a_1a_2+2a_0+a_2}{2a_1+a_1^2a_2}.$$ Now one may assume the existence of $h(x)\in k[x]$, $\deg h\le g$, such that $h (2a_1+a_1^2a_2)=2a_0a_1a_2+2a_0+a_2$. Therefore, one obtains $$\displaystyle a_0=\frac{h(2a_1+a_1^2a_2)-a_2}{2(a_2a_1+1)}=\frac{ha_1}{2}+\frac{ha_1-a_2}{2(a_1a_2+1)}.$$ Following the same ideas of the proof of Theorem \ref{thm1}, one may assume that $h=qa_2+r$, where $0\le \deg q\le g-1$ and $0\le \deg r< \deg a_2$. 
Setting $a_2=r a_1 -q$, one gets the equation
\[y^2=\left(rqa_1^2-q^2 a_1+ra_1 +q\right)^2+4(rqa_1 -q^2 +r ).\]
Choosing $r,q$ to be $\Q$-rationals and $q^2\ne r$ when $g$ is odd; and $q$ to be $\Q$-rational whereas $r$ to be a linear polynomial if $g$ is even, the order of the torsion divisor at infinity becomes $5 (g+1)/2$ if $g$ is odd; and $5g/2+2$ if $g$ is even, see Proposition \ref{prop:1}.

\begin{Example}
The following equation
\[C_t:\,y^2=(2(x^2+t)^2+(x^2+t)+1)^2+8(x^2+t)+4\] describes a family of hyperelliptic curves of genus $3$
over $\Q(t)$. The discriminant of the latter equation is given by $2^{48}\cdot 5^4\cdot (1 + t) (5 + 5 t + 4 t^3)$. If $t\in\Q\setminus\{-1\}$, then the class of the divisor \( (1 : 2 : 0) - (1 : -2 : 0) \) on the Jacobian of \( C_t \) is of order $10$.
\end{Example}
 
 \section{Explicit families}
In this section we provide explicit parametric families of hyperelliptic curves of genus $g\ge 3$ whose Jacobian possess $\Q$-rational torsion points of orders that does not appear in the literature. 
\begin{Theorem}
\label{thm13}
There exists infinitely many non-isomorphic hyperelliptic curves of genus $3$ over $\Q$ whose Jacobians posses a rational point of order $13$.
\end{Theorem}
\begin{Proof}
In Theorem \ref{thm1}, pick $r(x)=ax+b$, $a\ne 0$, to be a linear polynomial in $\Q[x]$, $u(x)=u\in\Q\setminus\{0\}$, and $a_1(x)=cx+d\in\Q[x]$, $c\ne 0$. When the discriminant of this curve is nonzero, we get a $5$-parameter family of hyperelliptic curves of genus $3$ satisfying the hypothesis of the theorem, where that divisor at infinity is $(1:a^2\cdot c^2\cdot u:0 )-(1:-a^2\cdot c^2\cdot u:0)$.  

In order to show that there are infinitely many such non-isomorphic curves; we set $a=1$, $b=t$, $c=d=u=1$. This provides a $1$-parameter subfamily $C^{13}_t$ of the aforementioned $5$-parameter family of curves. In addition, the discriminant of the defining equation is $-2^{28} h(t)$, where $$h(t)=36383 + 61280 t - 5762 t^2 - 15163 t^3 + 1662 t^4 + 2470 t^5 - 
  936 t^6 + 49 t^7 + 16 t^8.$$ In particular, if $t\in\Q$, then $C^{13}_t$ is a hyperelliptic curve. If two such curves $C^{13}_{t_1}$ and $C^{13}_{t_2}$, $t_1,t_2\in\Q$, are isomorphic over $\Q$, then $\Delta_{t_1}/\Delta_{t_2}$ is a $56$-th power of a non-zero rational number, \cite[\S 2]{Liu}, where $\Delta_{t_i}$ is the discriminant of the defining equation of $C^{13}_{t_i}$, $i=1,2$. A classical application of the Sieve of Eratosthenes can be used to prove that the integers $n$ for which $h(n)$ is $8$-th power free have positive density, \cite{Erdos}.   This implies that there are infinitely many non-isomorphic hyperelliptic curves in the family $C^{13}_t$.
\end{Proof}
In a similar fashion, one may prove the following theorem. 
\begin{Theorem}
\label{thm15}
There exists infinitely many non-isomorphic hyperelliptic curves of genus $4$ over $\Q$ whose Jacobians posses a rational point of order $15$.
\end{Theorem}
\begin{Proof}
In Theorem \ref{thm1}, we set $r(x)=a_1(x)=x+1$, $u(x)=x+t$. This gives rise to a $1$-parameter family $C_t^{15}$ of hyperelliptic curves of genus $4$ with a torsion divisor at infinity of order $15$. In addition, the discriminant of the defining equation is $2^{36}h(t)$ where {\footnotesize\begin{eqnarray*}
h(t)&=&25815184 + 4208616 t - 12580992 t^2 + 15705528 t^3 - 6568623 t^4 - 
 4422996 t^5 + 6829050 t^6 - 3728520 t^7\\ &+& 1148643 t^8 - 247912 t^9 + 
 48114 t^{10} - 8748 t^{11} + 729 t^{12}.
 \end{eqnarray*}}
 Now we conclude as in the proof of Theorem \ref{thm13}.
\end{Proof}

\begin{Remark}
\label{remm}
In the proof of Theorem \ref{thm13}, the equation describing the curve $C^{13}_{t}$ is of the form $(x+t)f_{t}(x)$ where the polynomial $f_{t}$ has generic Galois group $\textbf{S}_7$. The latter statement can be justified by seeing that the specialization $f_{1}$ has Galois group $\textbf{S}_7$. Hilbert's Irreducibility Theorem implies that the Galois groups of $f_{1}$ and $f_{t}$ over $\Q$ coincide for infinitely many rational values of $t$. According to Theorem \ref{thm1}, this means that the Jacobian $J_{C^{13}_{t}}$ of $C^{13}_{t}$ satisfies that $\operatorname{End}(J_{C^{13}_{t}}) = \Z$ and that $J_{C^{13}_{t}}$ is an absolutely simple abelian variety for infinitely many rational values of $t$. 

The same argument yields that the Jacobian of $C^{15}_{t}$ in the proof of Theorem \ref{thm15} is an absolutely simple variety for infinitely many rational values of $t$. 
\end{Remark}
We notice that for $g\ge 5$, there are at least two integers $N$ in the interval $[3g,4g+1]$ when $g$ is odd, $[3g,4g-1]$ when $g$ is even, satisfying the hypotheses of Corollary \ref{cor1}. 
\begin{Example}
Consider the following families of hyperelliptic curves of genus $5$
\begin{eqnarray*}
C^{17}_{s,t}&:&y^2 = 4 (t + x + x (t + x)^2 ( x^2+s)) + ( 
    x+t)^2 (x + ( x^2+s) (1 + x^2 ( x+t)))^2\\
C^{18}_{s,t,u}&:&y^2= ( x^2+s)^2 (t + u + x + u ( x+t)^2 ( x^2+s))^2 + 
 4 (s + x^2 + u ( x+t) ( x^2+s)^2) \\
C^{21}_{s,t,u}&:& y^2=4 (s + x + u ( x+s)^2 ( x^2+t)) + ( x+s)^2 (t + u + x^2 + 
    u ( x+s) (x^2+t)^2)^2
\end{eqnarray*}
corresponding to the triple $(\alpha,\beta,\gamma)=(1,1,2), (2,1,0),(1,2,0)$, respectively, in Theorem \ref{thm1}. In fact, for the triple $(1,1,2)$, we chose $a_1(x)=x$, $r(x)=x+t$, $u(x)=x^2+s$, which implies that the divisor at infinity is torsion of order $17$. 

For the triple $(1,2,0)$, we chose $a_1(x)=x+t$, $r(x)=x^2+s$, $u(x)=u\in \Q\setminus\{0\}$, yielding that the divisor at infinity is torsion of order $18$.
Finally, for the triple $(2,1,0)$, we chose $a_1(x)=x^2+t$, $r(x)=x+s$, $u(x)=u\in \Q\setminus\{0\}$, with the divisor at infinity being torsion of order $21$.

The same argument as in Remark \ref{remm} shows that the Jacobians of the curves $C_{s,t}^{17}$ and $C^{21}_{s,t,u}$ are absolutely simple varieties for infinitely many choices of $s,t,u$. 

\end{Example}

\section*{Conflict of Interest statement} On behalf of all authors, the corresponding author states that there is no conflict of interest.
\section*{Data availability statement} The authors declare that the data supporting the findings of this study are available within the paper.

\end{document}